\documentclass[twoside,11pt]{article}
\usepackage{graphicx,amsmath,amsthm,amssymb,latexsym,amsfonts,color}
\usepackage[bookmarksnumbered, colorlinks, plainpages]{hyperref}
\newtheorem{thm}{\bf Theorem}

\newtheorem{example}{\bf Example}

\numberwithin{equation}{section}
\def\h{\hspace{-0.2cm}}
\def\ao{\alpha_{\rm opt}}
\begin{document}
\title{Two-step scale-splitting method  for solving complex symmetric system of linear equations}
\author{{Davod Khojasteh Salkuyeh\footnote{Corresponding author}}\\[2mm]
\textit{{\small Faculty of Mathematical Sciences, University of Guilan, Rasht, Iran}} \\
\textit{{\small email: khojasteh@guilan.ac.ir}\textit{}}}
\date{}
\maketitle
{\bf Abstract.}
Based on the Scale-Splitting (SCSP) iteration method presented by Hezari et al. in (A new iterative method
for solving a class of complex symmetric system linear of equations, Numerical Algorithms 73 (2016) 927-955),
we present a new two-step iteration method, called TSCSP, for solving the complex symmetric system of linear equations $(W+iT)x=b$, where $W$ and $T$ are  symmetric positive definite and symmetric positive semidefinite matrices, respectively. It is shown that if the matrices $W$ and $T$ are symmetric positive definite, then the  method is unconditionally convergent. The optimal value of the parameter, which minimizes the spectral radius of the iteration matrix is also computed. Numerical {comparisons} of the TSCSP iteration method with the SCSP, the MHSS, the PMHSS and the GSOR methods
are given to illustrate the effectiveness of the method.
 \\[-3mm]


\pagestyle{myheadings}\markboth{~\hrulefill~D.K. Salkuyeh~}{Generalized SCSP iteration method~\hrulefill~}

\thispagestyle{empty}

\section{Introduction}\label{SEC1}
We consider the system of linear equations of the form
\begin{equation}\label{Eq1.01}
Az=(W+iT)z=b,\quad A\in \mathbb{C}^{n\times n},
\end{equation}
where  $z=x+iy$ and $b=f+ig$, in which the vectors $x,y,f$ and $g$ are in $\mathbb{R}^{n}$ and $i=\sqrt{-1}$. Assume that the matrices $W$ and $T$ are symmetric positive semidefinite  such that at least one of them, e.g., $W$, being positive definite. {This kinds of systems appear} in many applications
{including FFT-based} solution of certain time-dependent PDEs \cite{Bertaccini}, diffuse optical tomography \cite{Arridge},
algebraic eigenvalue problems \cite{Moro, Schmitt}, molecular scattering \cite{Poirier}, structural dynamics \cite{Feriani} and lattice
quantum chromodynamics \cite{Frommer}.

Several iteration methods have been presented to solve (\ref{Eq1.01}) in the literature. In \cite{HSS}, Bai et al.
proposed the Hermitian and skew-Hermitian splitting (HSS) method to solve the system of linear equations with non-Hermitian {positive definite coefficient} matrices. Next, Bai et al. presented the modified HSS iteration method (MHSS) for solving Eq (\ref{Eq1.01}) (see \cite{MHSS}). It is known that the matrix $A$ possesses the Hermitian/skew-Hermitian (HS) splitting
\[
A=H+S,
\]
where
\[
H=\frac{1}{2}(A+A^H)=W \quad \textrm{and} \quad S=\frac{1}{2}(A-A^H)=iT,
\]
in which  $A^H$ stands for the conjugate transpose of $A$. In this case, the MHSS iteration  method for solving \eqref{Eq1.01} can be written as follows.

\medskip

\noindent {\bf The MHSS iteration method}: \verb"Let" $z^{(0)} \in  {\mathbb{C }^{n}}$ \verb"be an initial guess". \verb"For" $k=0,1,2,\ldots$, until $\{z^{(k)}\}$ \verb"converges, compute" ${z^{(k+1)}}$  \verb"according to the following sequence":
\begin{equation*}
\begin{cases}
(\alpha I+W){z^{(k+\frac{1}{2})}}=(\alpha I-iT){z^{(k)}}+b, \\
(\alpha I+T){z^{(k+1)}}=(\alpha I+iW){z^{(k+\frac{1}{2})}}-ib,
\end{cases}
\end{equation*}
\verb"where" $\alpha$ \verb"is a given positive constant and" $I$ \verb"is the identity matrix".

\medskip

In \cite{MHSS}, it has been shown that when the matrices $W$ and $T$ are symmetric positive definite and symmetric positive semidefinite,
respectively, then the MHSS iteration method is convergent.  Since both of the matrices  $\alpha I+T$ and $\alpha I+W$ are symmetric positive
definite, the two linear subsystems involved in each step can be solved exactly by the Cholesky factorization \cite{Saadbook} of the coefficient
matrices or inexactly by the conjugate gradient (CG) method \cite{CG}. A preconditioned version of the MHSS (called PMHSS) method has been presented by Bai et
al. in \cite{PMHSS} which can be written as following.

\medskip

\noindent {\bf The PMHSS iteration method}: \verb"Let" $z^{(0)} \in  {\mathbb{C }^{n}}$ \verb"be an initial guess". \verb"For" $k=0,1,2,\ldots$, until $\{z^{(k)}\}$ \verb"converges, compute" ${z^{(k+1)}}$  \verb"according to the following sequence":
\begin{equation*}
\begin{cases}
(\alpha V+W){z^{(k+\frac{1}{2})}}=(\alpha V-iT){z^{(k)}}+b, \\
(\alpha V+T){z^{(k+1)}}=(\alpha V+iW){z^{(k+\frac{1}{2})}}-ib,
\end{cases}
\end{equation*}
\verb"where" $\alpha$ \verb"is a given positive constant and" $V$ \verb"is a symmetric positive definite".

\medskip

{If the matrix $V$ is specified to be the identity matrix $I$, then PMHSS is simplified as MHSS. Theoretical analysis in \cite{PMHSS} has shown that the MHSS iteration converges
to the unique solution of the complex symmetric linear system \eqref{Eq1.01} for any initial guess.}
In practice the matrix $V$ is set to $W$. Numerical results presented in \cite{PMHSS} show that the PMHSS iteration method outperforms MHSS. The subsystems appeared in the PMHSS method can be treated as the methods described for the MHSS method. Several variants of the HSS method have been presented by Wu in \cite{WuNLA}.

In \cite{GSOR}, Salkuyeh et al. proposed the GSOR iteration method to solve Eq. (\ref{Eq1.01}) which can be
written as following.

\medskip

\noindent {\bf The GSOR iteration method}: \verb"Let" $(x^{(0)}; y^{(0)}) \in  {\mathbb{R }^{n}}$ \verb"be an initial guess. For" $k=0,1,2,\ldots$, \verb"until" $\{(x^{(k)};y^{(k)})\}$ \verb"converges, compute" ${(x^{(k+1)};y^{(k+1)})}$  \verb"according to the" \\ \verb"following sequence"
\begin{equation*}
\begin{cases}
W x^{(k+1)}=(1-\alpha)Wx^{(k)}+\alpha T y^{(k)}+\alpha f, \\
Wy^{(k+1)}=-\alpha T x^{(k+1)}+(1-\alpha)W y^{(k)}+ \alpha g,
\end{cases}
\end{equation*}
\verb"where" $\alpha$ \verb"is a given positive constant".

\medskip

In \cite{GSOR} it has been shown that if $W$ is symmetric positive definite and $T$ is symmetric then the GSOR method is convergent if and only
if
\[
0<\alpha<\frac{2}{1+\rho(W^{-1}T)}.
\]
The optimal value of the parameter $\alpha$ in the GSOR method was also obtained in \cite{GSOR}.
In each iterate of the GSOR method, two subsystems with the coefficient matrix $W$ should be solved. {Since $W$ is symmetric positive definite, the Cholesky factorization of $W$ or the CG iteration  method can be utilized for solving these systems.}  In contrast with the MHSS method, GSOR is a method based on real arithmetic. Numerical results of the GSOR method showed that  in general, it outperforms the MHSS method \cite{GSOR}.

{Recently, Hezari et al. in \cite{SCSP} have presented a new iteration method called  Scale-Splitting (SCSP) for solving (\ref{Eq1.01}), which serves the SCSP preconditioner.  In \cite{SCSP} it was shown that the application of the GSOR method to the preconditioned system in conjunction with the SCSP preconditioner is very efficient. In this paper, we present a two-step SCSP method (TSCSP) for solving \eqref{Eq1.01} and compare
it with the SCSP, the MHSS, the PMHSS and the GSOR methods.}

Throughout the paper, for a square matrix $A$,  $\rho(A)$ and $\sigma(A)$ stand for the spectral radius and spectrum of $A$, respectively.

{This paper is organized as follows. Section \ref{SEC2} describes the TSCSP iteration method. Convergence of the TSCSP method is investigated in Section \ref{SEC3}. Numerical experiments are given in Section \ref{SEC4}. Some concluding remarks are presented in Section \ref{SEC5}.
}

\section{The TSCSP iteration method}\label{SEC2}

Let $\alpha>0$. As the SCSP method, we multiply both sides of Eq. (\ref{Eq1.01}) by $\alpha-i$ to get the equivalent system
{\begin{equation}\label{Eq2.01}
(\alpha-i)Az=(\alpha-i)b,
\end{equation}}
where $i=\sqrt{-1}$. Then, we split the coefficient matrix of the system (\ref{Eq2.01}) as
\[
(\alpha-i)A=(\alpha W+T)-i(W-\alpha T).
\]
Using this splitting, we rewrite system (\ref{Eq2.01}) as the fixed-point equation
\begin{equation}\label{Eq2.02}
(\alpha W+ T)z=i(W-\alpha T)z+(\alpha -i)b.
\end{equation}
On the other hand, we multiply both sides of Eq. (\ref{Eq1.01}) by $1-\alpha i$ to obtain the equivalent system
\[
(1-\alpha i)Ax=(1-\alpha i)b,
\]
which yields the fixed point equation
\begin{equation}\label{Eq2.04}
(W+\alpha T)z=i(\alpha W- T)z+(1-\alpha i)b.
\end{equation}
From Eqs. (\ref{Eq2.02}) and  (\ref{Eq2.04}),  we now state the TSCSP algorithm as follows.

\bigskip

\noindent {\bf The TSCSP iteration method:} \verb"Given an initial guess" $z^{(0)}$, \verb"for" $k=0, 1, 2,\ldots$, \verb"until" $z^{(k)}$ \verb"converges, compute"
\begin{eqnarray*}
(\alpha W+ T)z^{(k+\frac{1}{2})} &\h=\h& i(W-\alpha T)z^{(k)}+(\alpha -i)b,\\
(W+\alpha T)z^{(k+1)}  &\h=\h& i(\alpha W- T)z^{(k+\frac{1}{2})}+(1-\alpha i)b,
\end{eqnarray*}
\verb"where" $\alpha>0$.

\bigskip

In each iterate of the TSCSP iteration method two subsystems with the coefficient matrices  $\alpha W+ T$ and $W+\alpha T$ should be solved.
Both of these matrices are symmetric positive definite. Hence, the subsystems can be solved directly by the Cholesky factorization
or the conjugate gradient method inexactly. {Obviously, an iteration step of the TSCSP method is completely equivalent to two iteration steps of the SCSP method in terms of the computation cost.}

Computing the vector $z^{(k+\frac{1}{2})}$ from the first step of the TSCSP method
and substituting it in the second step gives the following stationary method
\begin{align}
\label{TSCSP-St}
z^{k+1}=\mathcal{G}_{\alpha} z^{(k)}+c_{\alpha},
\end{align}
where
\begin{eqnarray*}\label{G}
\mathcal{G}_{\alpha}&\h=\h& (W+\alpha T)^{-1}(T-\alpha W)(\alpha W+T)^{-1}(W-\alpha T),
\end{eqnarray*}
and
\[
c_{\alpha}=2\alpha(W+\alpha T)^{-1} (W-iT) (\alpha W+T)^{-1} b.
\]
Since $W$ is Hermitian positive definite, we can write $W=W^{\frac{1}{2}}W^{\frac{1}{2}}$, with  $W^{\frac{1}{2}}$ being
Hermitian positive definite. Using this fact, the iteration {matrix $\mathcal{G}_{\alpha}$ can} be written as
\[
\mathcal{G}_{\alpha}=W^{-\frac{1}{2}} (I+\alpha S)^{-1}(S-\alpha I)(\alpha I+S)^{-1}(I-\alpha S)  W^{\frac{1}{2}},
\]
where $S=W^{-\frac{1}{2}} T W^{-\frac{1}{2}}$. Obviously, $S$ is Hermitian positive semidefinite and as a result its eigenvalues are nonnegative. Letting
\begin{equation}\label{Ghat}
\hat{\mathcal{G}}_{\alpha}=(I+\alpha S)^{-1}(S-\alpha I)(\alpha I+S)^{-1}(I-\alpha S),
\end{equation}
it follows that
\begin{equation}\label{GGhatSim}
\mathcal{G}_{\alpha}=W^{-\frac{1}{2}}  \hat{\mathcal{G}}_{\alpha} W^{\frac{1}{2}}.
 \end{equation}
This shows that the matrices $\mathcal{G}_{\alpha}$ and $\hat{\mathcal{G}}_{\alpha}$ are similar. It is worth noting that, to prove the convergence of the TSCSP iteration method,  working with the matrix $\hat{\mathcal{G}}_{\alpha}$ would be easier than the matrix $\mathcal{G}_{\alpha}$.

Letting
\begin{eqnarray*}
M_{\alpha}& \h=\h & \frac{1}{2\alpha}(T+\alpha W)(W-iT)^{-1}(W+\alpha T),\\
N_{\alpha}& \h=\h & \frac{1}{2\alpha}(T-\alpha W)(W-iT)^{-1}(W-\alpha T),
\end{eqnarray*}
it follows that
\[
A=M_{\alpha}-N_{\alpha} \quad \textrm{and} \quad \mathcal{G}_{\alpha}=M_{\alpha}^{-1}N_{\alpha}.
\]
Hence, the TSCSP iteration method is induced by the matrix splitting $A=M_{\alpha}-N_{\alpha}$.  It follows that the matrix $M_{\alpha}$
can be used as a preconditioner for the system (\ref{Eq1.01}), which is referred to as the TSCSP preconditioner.

\section{Covergence of the TSCSP iteration method}\label{SEC3}

For the convergence of the TSCSP iteration method, all we need to do is to provide conditions under which $\rho(\mathcal{G}_{\alpha})<1$.
To do this, we state and prove the following theorem.
\begin{thm}\label{Thm1}
Assume that the matrices $W\in \mathbb{R}^{n\times n}$ and $T\in \mathbb{R}^{n\times n}$ are symmetric positive definite. Then, for every $\alpha>0$,
$\rho(\mathcal{G}_{\alpha})<1$. That is, the TSCSP iteration method is convergent {for every} $\alpha>0$.
\begin{proof}
Since both of the matrices $W$ and $T$ are symmetric positive definite, we conclude that the eigenvalues of $S=W^{-\frac{1}{2}} T W^{-\frac{1}{2}}$ are positive. Therefore, from Eqs. (\ref{Ghat}) and (\ref{GGhatSim})  we obtain
\begin{eqnarray*}
\rho(\mathcal{G}_{\alpha}) &\h=\h& \rho(\hat{\mathcal{G}}_{\alpha})\\
                           &\h=\h& \rho\left((I+\alpha S)^{-1}(S-\alpha I)(\alpha I+S)^{-1}(I-\alpha S)\right)\\
                           &\h=\h& \max_{\mu_j\in\sigma(S)}\left|\frac{(\mu_j-\alpha)(1-\alpha \mu_j)}{(\mu_j+\alpha)(1+\alpha \mu_j)}\right|\\
                           &\h=\h& \max_{\mu_j\in\sigma(S)}\left|\frac{\mu_j-\alpha}{\mu_j+\alpha}\right|\left|\frac{1-\alpha \mu_j}{1+\alpha \mu_j}\right|\\
                           &\h<\h&1,
\end{eqnarray*}
which completes the proof.
\end{proof}
\end{thm}

 The next theorem presents the optimal value of the parameter $\alpha$ in the TSCSP iteration method.

\begin{thm}\label{OptThm}
Let $W\in \mathbb{R}^{n\times n}$ and $T\in \mathbb{R}^{n\times n}$ be symmetric positive definite matrices. Let also $\mu_i$, $i=1,2,\ldots,n$ be the eigenvalues of $S=W^{-\frac{1}{2}} T W^{-\frac{1}{2}}$.\\
(a) If $0<\mu_1\leq\mu_2\leq \cdots \leq \mu_n\leq 1$ or $1\leq \mu_1\leq\mu_2\leq \cdots \leq \mu_n$  set $\gamma=\mu_1$ and $\delta=\mu_n$.\\
(b) If $0<\mu_1<\cdots \leq \mu_k \leq 1 \leq \mu_{k+1}\leq \cdots \leq \mu_n$, then if $\mu_1\mu_n\geq 1$, set $\gamma=\mu_{k+1}$ and $\delta=\mu_n$, otherwise $\gamma=\mu_{1}$ and $\delta=\mu_k$.\\
Then, the optimal values of $\alpha$ in the TSCSP iteration method are given by
\[
\alpha_{opt}^{\pm}=\rm{arg}\min_{\hspace{-0.2cm}\alpha>0}\rho(\mathcal{G}_{\alpha})=\frac{1}{2}\left(\eta\pm\sqrt{\eta^2-4}\right),
\]
where
\[
\eta=\sqrt{\frac{(1+\gamma^2)(1+\delta^2)}{\gamma\delta}}.
\]
Moreover
\[
\rho\left(\mathcal{G}_{\alpha_{opt}^{\pm}}\right)=\left|\frac{\delta^2-\eta\delta+1}{\delta^2+\eta\delta+1}\right|.
\]
\begin{proof}
We prove the theorem for the case  $0<\mu_1\leq\mu_2\leq \cdots \leq \mu_n\leq 1$. Other cases can be similarly proved. Having Theorem \ref{Thm1} in mind, we define the function $\lambda_{\mu}(\alpha)$ as
\[
\lambda_{\mu}(\alpha)=\frac{(\mu-\alpha)(1-\alpha \mu)}{(\mu+\alpha)(1+\alpha \mu)}=\frac{(\alpha-\mu)(\alpha-\frac{1}{\mu})}{(\alpha+\mu)(\alpha+\frac{1}{\mu})}=\frac{\alpha^2-\beta\alpha+1}{\alpha^2+\beta\alpha+1},
\]
where $\beta=\mu+1/\mu$. Obviously, $\beta\geq 2$. Some properties of the function $\lambda_{\mu}(\alpha)$ are given as follows.
This function passes through the points $(0,1)$, $(\mu,0)$ and $(1/\mu,0)$ and $y=1$ is its vertical asymptote.
By direct computations, we get
\[
\frac{d}{d\alpha} \lambda_{\mu}(\alpha)=\frac{2\beta(\alpha^2-1)}{(\alpha^2+\beta\alpha+1)^2}.
\]
Therefore, the only critical point of the function $\lambda_{\mu}(\alpha)$ in the {interval $(0,+\infty)$} is $\alpha=1$.
Moreover, the function $\lambda_{\mu}(\alpha)$ in the interval $(0,1)$ is strictly decreasing and in the interval $(1,+\infty)$ strictly increasing.
On the other hand, it is easy to see that
\[
\frac{d^2}{d\alpha^2}\lambda_{\mu}(\alpha)=\frac{4\beta(3\alpha-\alpha^3+\beta)}{(\alpha^2+\beta\alpha+1)^3}.
\]
It is not difficult to see that there exists $\alpha_0>1$ such that
\[
\frac{d^2}{d\alpha^2}\lambda_{\mu}(\alpha)\geq 0,\quad  \textrm{for}\quad {0<\alpha<\alpha_0},\quad \textrm{and} \quad \frac{d^2}{d\alpha^2}\lambda_{\mu}(\alpha)\leq 0,\quad  \textrm{for}\quad \alpha\geq \alpha_0.
\]
Hence, the function is concave upward in the interval $(0,\alpha_0)$ and concave downward in the interval $(\alpha_0,+\infty)$.
Using the above information about the function $\lambda_{\mu}(\alpha)$, for a given $\mu$ the function has been displayed in  Figure \ref{Fig1}.
\begin{figure}
\centerline{\hspace{2cm}\includegraphics[height=9cm,width=14cm]{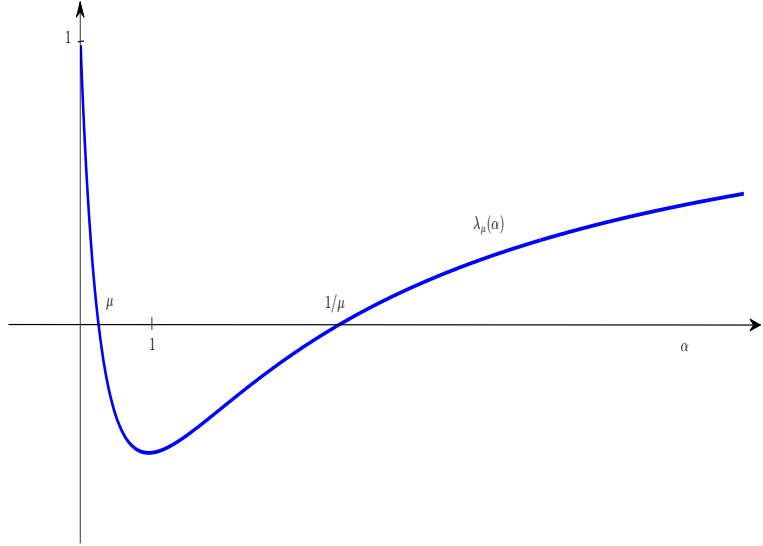}}\vspace{-1.5cm}
\caption{Graph of $\lambda_{\mu}(\alpha)$ for a given $\mu<1$ and $\alpha>0$.}\label{Fig1}
\end{figure}
For $0<\mu_1\leq\mu_2\leq \mu_3\leq 1$ the absolute value of the functions $\lambda_{\mu_i}(\alpha)$, $i=1,2,3$, have been displayed in Figure \ref{Fig2}.
\begin{figure}
\centerline{\hspace{1cm}\includegraphics[height=9cm,width=15cm]{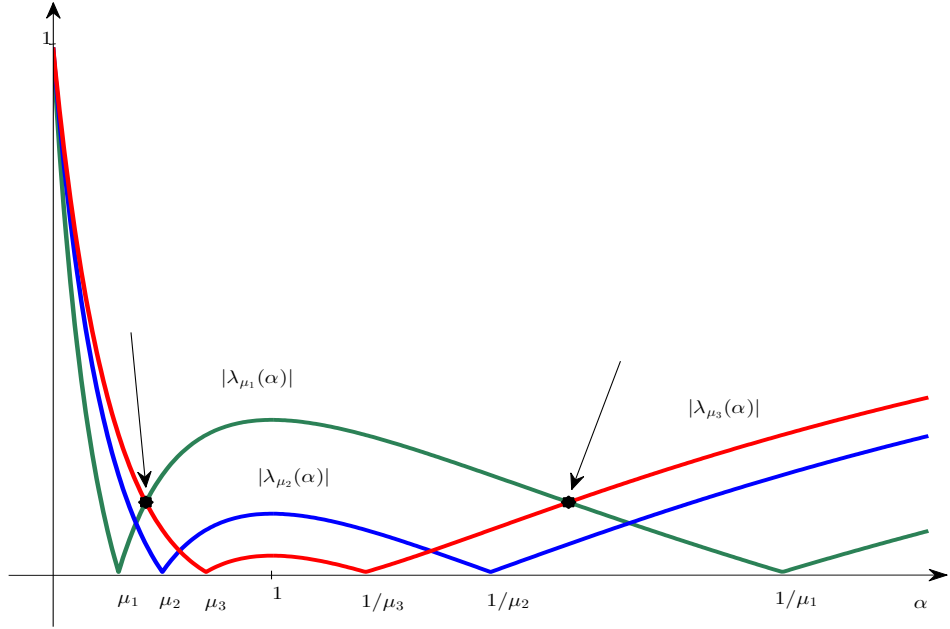}}
\caption{Graph of $|\lambda_{\mu}(\alpha)|$ for $\mu=\mu_1,\mu_2,\mu_3$  where $\mu_1 < \mu_2 < \mu_3 \leq 1$ and for $\alpha>0$.}\label{Fig2}
\end{figure}
The optimal points have been pointed out in the figure by two bullets. As seen the optimal values of $\alpha$ are obtained by intersecting the functions $|\lambda_{\mu_1}(\alpha)|$ and $|\lambda_{\mu_3}(\alpha)|$. Hence, in the general case, if we set $\gamma=\mu_1$ and $\delta=\mu_n$, then the optimal value of $\alpha$ satisfies the relation
\[
\frac{(\gamma-\alpha_{opt})(1-\ao\gamma)}{(\gamma+\ao)(1+\ao\gamma)}
=-\frac{(\delta-\ao)(1-\ao\delta)}{(\delta+\ao)(1+\ao\delta)},
\]
which is equivalent to
\[
\frac{-\gamma^2+\eta\gamma-1}{\gamma^2+\eta\gamma+1}=\frac{\delta^2-\eta\delta+1}{\delta^2+\eta\delta+1},
\]
where $\eta=\ao+1/\ao$. Since $\ao>0$, if we solve the latter equation for $\eta$, then we deduce
\begin{equation}\label{Eta}
\ao+\frac{1}{\ao}=\eta=\sqrt{\frac{(1+\gamma^2)(1+\delta^2)}{\gamma\delta}}.
\end{equation}
Therefore, to compute $\ao$ we need to solve the quadratic equation
\begin{equation}\label{QuadEq}
\ao^2-\eta \ao+1=0.
\end{equation}
Since $1+\gamma^2\geq 2 \gamma$ and $1+\delta^2\geq 2\delta$, it follows from Eq. (\ref{Eta}) that $\eta\geq 2$ and as a result the discriminant
of Eq. (\ref{QuadEq}) is nonnegative, i.e., $\Delta=\eta^2-4\geq 0$.
By solving  Eq. (\ref{QuadEq}) the desired result is obtained.
\end{proof}
\end{thm}

Some comments can be posed here. For the TSCSP iteration method, it follows from Theorem \ref{OptThm} that, for each problem, there are in general two optimal values for the parameter $\alpha$, one of them is less than or equal to $1$ and the other one is greater than or equal to 1.
{This is an interesting property of the TSCSP iteration method, since the optimal value of the method can be sought in the interval $(0,1]$.} As Theorem \ref{OptThm} shows, the optimal value of  $\alpha$ minimizes the spectral radius of the iteration matrix of the TSCSP iteration method, whereas in the HSS iteration method (as well as the MHSS iteration method) it minimizes an upper bound of the spectral radius of the iteration matrix.
When $W$ is symmetric positive definite and $T$ is symmetric positive semidefinite, the HSS, the MHSS and the GSOR methods are convergent. However, according to Theorem \ref{Thm1} for the convergence of the TSCSP iteration method both of the matrices $W$ and $T$ should be symmetric positive definite. In fact, when $T$ is singular, at least one of the eigenvalues of $T$ is zero and therefore it follows that $\rho\left(\mathcal{G}_{\alpha}\right)=1$ and the method fails to converge.

\section{Numerical experiments}\label{SEC4}

In this section we present four examples  to demonstrate the feasibility and effectiveness
of the TSCSP iteration method.  Three of them have been chosen from  \cite{MHSS} and one of them has been made artificially.   We compare the numerical results of the TSCSP iteration method with those of the SCSP, the MHSS, the PMHSS and the GSOR
methods from the point of view of both the number of iterations (denoted by ``Iter") and the total computing time (in seconds, denoted by ``CPU").
In each iteration of these methods, we use the Cholesky factorization of the coefficient matrices to solve the sub-systems.
The reported CPU times  are the sum of the CPU time for the convergence of the method and the CPU time for computing the Cholesky factorization.
It is necessary to mention that for solving the symmetric positive definite system of linear equations we have used the sparse Cholesky
factorization incorporated with the symmetric approximate minimum degree reordering \cite{Saadbook}. To do so, we have used the \verb"symamd.m"
command of \textsc{Matlab} Version 7.

All the numerical experiments  were computed in double precision using some  \textsc{Matlab}
 codes  on a Pentium 4 Laptop, with a 1.80 GHz CPU
and 2.40GB of RAM. We use a null vector as an initial guess and the  stopping criterion
\[
\frac{\|b-Az^{(k)}\|_{2}}{\| b\|_{2}}<10^{-6},
\]
is always used, where $z^{(k)}$ is the computed solution at iterate $k$.


\begin{example} \label{ex1}\rm  (See \cite{MHSS})
We consider the system (\ref{Eq1.01}) with
\[
W=K+\frac{3-\sqrt{3}}{\tau}I,\quad T=K+\frac{3+\sqrt{3}}{\tau}I,
\]
where $\tau$ is the time step-size and $K$ is the five-point centered difference matrix approximating the negative Laplacian
operator $L\equiv -\Delta$ with homogeneous Dirichlet boundary conditions, on a uniform mesh in the unit square
$[0, 1]\times[0,1]$ with the mesh-size $h=1/(m+1)$. The matrix $K\in\mathbb{R}^{n\times n}$ possesses the tensor-product form
$K=I\otimes V_{m}+V_{m}\otimes I$, with $V_{m}=h^{-2}{\rm tridiag}(-1,2,-1)\in \mathbb{R}^{m\times m}$. Hence, $K$ is an
${n\times n}$ block-tridiagonal matrix, with $n=m^{2}$. The right-hand side vector $b$ with its $j$th entry $b_{j}$ being given by
\[
b_{j}=\frac{(1-i)j}{\tau(j+1)^{2}},\quad  j=1,2,\ldots,n.
\]
In our tests, we take $\tau=h$. Furthermore, we normalize coefficient matrix and right-hand side by multiplying both by $h^{2}$.
\end{example}


\begin{example}\label{ex2} \rm (See \cite{MHSS})
We consider the system $Az=(W+iT)z=b$, with
$$W=-\omega^{2}M+K,\quad T=\omega C_{V}+C_{H},$$
 where $M$ and $K$ are the inertia and the stiffness matrices, $C_{V}$ and $C_{H}$ are the
 viscous and the hysteretic damping matrices, respectively, and $\omega$ is the driving circular frequency.
 We take $C_{H}=\mu K$ with $\mu$ a damping coefficient, $M=I$, $C_{V}=10I$, and $K$ the five-point centered
 difference matrix approximating the negative Laplacian operator with homogeneous Dirichlet boundary conditions,
 on a uniform mesh in the unit square $[0, 1]\times[0, 1]$ with the mesh-size $h=1/(m+1)$.
 The matrix $K\in \mathbb{R}^{n\times n}$ possesses the tensor-product form $K=I\otimes V_{m}+V_{m}\otimes I$,
 with $V_{m}=h^{-2}{\rm tridiag}(-1,2,-1)\in \mathbb{R}^{m\times m}$. Hence, $K$ is an ${n\times n}$ block-tridiagonal matrix,
 with $n=m^{2}$. In addition, {we set $\omega=4$}, $\mu=0.02$, and the right-hand side vector $b$ to be $b=(1 + i)A{\bf 1}$, with ${\bf 1}$
 being the vector of all entries equal to $1$. As before, we normalize the system by multiplying both sides
through by $h^{2}$.
\end{example}


\begin{example}\label{ex3}\rm (See \cite{MHSS})
Consider the system of linear equations $(W+iT)x=b$, with
\[
T=I\otimes V+V\otimes I \quad {\rm and} \quad W=10(I\otimes V_{c}+V_{c}\otimes I)+9(e_{1}e_{m}^{T}+e_{m}e_{1}^{T})\otimes I,
\]
where $V={\rm tridiag}(-1,2,-1)\in \mathbb{R}^{m\times m}$, $V_{c}=V-e_{1}e_{m}^{T}-e_{m}e_{1}^{T}\in \mathbb{R}^{m\times m}$
and $e_{1}$  and $e_{m}$  are the first and last unit vectors in $\mathbb{R}^{m}$, respectively. We take the right-hand side
vector $b$ to be $b=(1 + i)A\textbf{1}$, with $\textbf{1}$ being the vector of all entries equal to $1$.

Here $T$ and $W$ correspond to the five-point centered difference matrices approximating the negative Laplacian operator with
homogeneous Dirichlet boundary conditions and periodic boundary conditions, respectively, on a uniform mesh in the unit square
$[0, 1]\times[0, 1]$ with the mesh-size $h=1/(m+1)$.
\end{example}

\begin{example}\label{ex4}\rm
In this example we artificially set
\[
W={\rm tridiag}(-1+\theta_1,2,-1+\theta_1)\in \mathbb{R}^{n\times n},\quad  T={\rm tridiag}(-1+\theta_2,2,-1+\theta_2)\in \mathbb{R}^{n\times n},
\]
where $\theta_1,\theta_2\in \mathbb{R}$. We solve the system $(W+iT)z=b$ by the methods, where $b=A\textbf{1}$ with $\textbf{1}$
being a vector of all ones. In our test problem we set $\theta_1=1.5$, $\theta_2=0.2$.
\end{example}


\begin{table}
\centering\caption{Numerical results for Example \ref{ex1}. \label{Tbl1}}
\bigskip
\begin{tabular}{lllllllll} \hline
Method            &$n=m^2$& $32^2$&$64^2$ &$128^2$&$256^2$&$512^2$& $1024^2$       \\ \hline \vspace{-0.3cm} \\ \vspace{0.0cm}

TSCSP             & $\ao$ & 0.46  & 0.46  & 0.46  & 0.46  & 0.46  & 0.46 \\
                  & Iter  & 7     & 7     & 7     &  7    & 7     & 7 \\
                  & CPU   & 0.02  & 0.04  & 0.21  & 0.96  & 4.88  & 24.94 \\[0.2cm]

SCSP              & $\ao$ & 0.65  & 0.65  & 0.65 & 0.65  & 0.65  & 0.65 \\
                  & Iter  & 9     & 9     & 9    & 9     & 9     & 9 \\
                  & CPU   & 0.01  & 0.02  & 0.11 & 0.50  & 2.87  & 15.17 \\[0.2cm]

MHSS             & $\ao$  & 0.78  & 0.55 & 0.40  & 0.30  & 0.21  & 0.15 \\
                  & Iter    & 53    & 72   & 98  & 133   & 181   & 249 \\
                  & CPU   & 0.03  & 0.16 & 1.45  & 10.60 & 82.13 & 632.78 \\[0.2cm]

PMHSS             & $\ao$ & 1.36  & 1.35 & 1.05  & 1.05  & 1.05  & 1.05 \\
                  & Iter  & 21    & 21   & 21    & 21    & 20    & 20 \\
                  & CPU   & 0.02  & 0.06 & 0.36  & 2.00  & 10.59 & 62.51 \\[0.2cm]

GSOR             & $\ao$  & 0.495 & 0.457& 0.432 & 0.418 & 0.412 & 0.411 \\
                  & Iter  &   22  & 24   & 26    & 27    & 27    & 27 \\
                  & CPU   & 0.01  & 0.06 & 0.43  & 2.42  & 13.69 & 66.73  \\ \hline
 \end{tabular}
\end{table}

Tables \ref{Tbl1}-\ref{Tbl4} present the numerical results for Examples \ref{ex1}-\ref{ex4} with different values of $n$.
In these tables, the optimal value of the parameter $\alpha$ (denoted by $\ao$) for the methods along {with the number of} the
iterations and the CPU time  for the convergence have been given. {For all the methods the optimal values of the
parameters were obtained experimentally}, except for the GSOR itration method, that the formula for the optimal value of $\alpha$ given in \cite{GSOR} was used.

Table \ref{Tbl1} contains the numerical results for Example \ref{ex1}. We see that the number of iterations of TSCSP is always less than those  of the other methods. However, from the CPU time point of view the TSCSP iteration method outperforms the MHSS, the PMHSS and the GSOR methods, but it can not compete with the SCSP method.  It is also seen that by increasing the dimension of the problem the number of iterations of the TSCSP and the SCSP iteration methods remain constant. It is worthwhile to note here that the number of iterations of the TSCSP, the SCSP,  the PMHSS and the GSOR methods is not too sensitive with respect to the size of the problems. However, the number of iterations of MHSS grows rapidly with the problem size.  Surprisingly,
we see from Table \ref{Tbl1} that for the TSCSP and the SCSP methods,  the optimal value of the parameter $\alpha$ are approximately 0.46 and 0.65, respectively.

Numerical results for Examples \ref{ex2}-\ref{ex4} are listed in Tables \ref{Tbl2}-\ref{Tbl4}, respectively.  As we observe, the TSCSP  outperforms the other methods from both the number of iterations and the CPU time point of view. Many of the comments and observations which we gave for Example \ref{ex1} can also be posed for Examples \ref{ex2}-\ref{ex4}.


\begin{table}[!ht]
\centering\caption{Numerical results for Example \ref{ex2}. \label{Tbl2}}
\bigskip
\begin{tabular}{lllllllll} \hline
Method            &$n=m^2$& $32^2$&$64^2$ &$128^2$&$256^2$&$512^2$& $1024^2$      \\ \hline \vspace{-0.3cm} \\ \vspace{0.0cm}

TSCSP             & $\ao$ & 0.11  & 0.09  & 0.08  & 0.07  & 0.07  & 0.06\\
                  & Iter  & 24    & 26    & 26    &  25   & 24    & 22  \\
                  & CPU   & 0.02  & 0.08  & 0.44  & 2.25  & 11.88 & 59.23\\[0.2cm]

SCSP              & $\ao$ & 1.07  & 1.09  & 1.10  & 1.10  & 1.11  & 1.12\\
                  & Iter  & 104   & 107   & 106   & 102   & 92    & 92\\
                  & CPU   & 0.03  & 0.13  & 0.80  & 4.08  & 22.13 & 120.35 \\[0.2cm]

MHSS             & $\ao$  & 0.08  & 0.04  & 0.02  & 0.01  & 0.005 & 0.003\\
                  & Iter  & 38    & 51    & 81    & 138   & 249   & 452\\
                  & CPU   & 0.03  & 0.13  & 1.19  & 10.77 & 109.31& 1063.81  \\[0.2cm]

PMHSS             & $\ao$ & 0.73  & 0.74  & 0.75  & 0.76  & 0.77  & 0.78\\
                  & Iter  & 36    & 38    & 38    & 38    & 38    & 38\\
                  & CPU   & 0.04  & 0.12  & 0.61  & 3.15  & 17.85 & 98.15 \\[0.2cm]

GSOR             & $\ao$  & 0.167 & 0.167 & 0.167 & 0.167 & 0.167 & 0.167\\
                  & Iter  & 76    & 76    & 76    & 76    & 76    & 76\\
                  & CPU   & 0.04  & 0.15  & 1.01  & 5.39  & 31.30 & 182.58\\ \hline
 \end{tabular}
\end{table}

\begin{table}[!ht]
\centering\caption{Numerical results for Example \ref{ex3}. \label{Tbl3}}
\bigskip
\begin{tabular}{lllllllll} \hline
Method            &$n=m^2$& $32^2$&$64^2$ &$128^2$&$256^2$&$512^2$& $1024^2$     \\ \hline \vspace{-0.3cm} \\ \vspace{0.0cm}

TSCSP             & $\ao$  & 0.23  & 0.23  & 0.23  & 0.23  &  0.16   & 0.11 \\
                  & Iter   & 13    & 13    & 13    &  13   &   16    & 23  \\
                  & CPU    & 0.02  & 0.07  & 0.40  & 2.19  &  15.35  & 151.39 \\[0.2cm]

SCSP              & $\ao$  & 1.92  & 1.44  & 1.15  & 1.02  & 0.96   & 0.93  \\
                  & Iter   & 15    & 25    & 40    & 59    & 78     & 94 \\
                  & CPU    & 0.01  & 0.07  & 0.50  & 4.07  & 33.06  & 233.34\\[0.2cm]

MHSS             & $\ao$   & 1.05  & 0.55  & 0.27  & 0.14  & 0.07   &  \\
                  & Iter   & 75    & 128   & 241   & 458   & 869    & $\dag$ \\
                  & CPU    & 0.04  & 0.36  & 4.49  & 47.64 & 528.54 &  \\[0.2cm]

PMHSS             & $\ao$  & 0.42  & 0.57  & 0.78  & 0.73  & 0.73   & 0.78 \\
                  & Iter   & 30    & 30    & 30    & 30    & 32     & 33 \\
                  & CPU    & 0.02  & 0.13  & 0.76  & 4.39  & 28.52  & 180.01  \\[0.2cm]

GSOR             & $\ao$   & 0.776 & 0.566 & 0.351 & 0.193 & 0.104  & 0.0545 \\
                  & Iter   & 11    & 20    & 33    & 64    & 129    & 261 \\
                  & CPU    & 0.01  & 0.09  & 0.72  & 8.24  & 112.23 & 1160.90 \\ \hline
 \end{tabular}
\end{table}


\begin{table}[!ht]
\centering\caption{Numerical results for Example \ref{ex4}. \label{Tbl4}}
\bigskip
\begin{tabular}{lllllllll} \hline
Method            &  $n$   & $32^2$  & $64^2$  & $128^2$ & $256^2$ & $512^2$ & $1024^2$            \\ \hline \vspace{-0.3cm} \\ \vspace{0.0cm}

TSCSP             & $\ao$  & 0.22  & 0.22  & 0.20  & 0.20  & 0.20  & 0.19 \\
                  & Iter   & 11    & 10    & 10    & 10    & 9     & 8    \\
                  & CPU    & 0.01  & 0.02  & 0.05  & 0.23  & 1.07  & 3.91 \\[0.2cm]

SCSP              & $\ao$  & 1.34  & 1.36  & 1.36  & 1.37  & 1.42  & 1.45 \\
                  & Iter   & 26    & 25    & 24    & 21    & 22    & 21 \\
                  & CPU    & 0.01  & 0.02  & 0.05  & 0.23  & 1.15  & 4.36 \\[0.2cm]

MHSS             & $\ao$   & 1.70  & 1.70  & 1.70  & 1.70  & 1.70  & 1.70 \\
                  & Iter   & 28    & 28    & 28    & 28    & 28    & 28 \\
                  & CPU    & 0.02  & 0.03  & 0.10  & 0.50  & 2.76  & 10.54 \\[0.2cm]

PMHSS             & $\ao$  & 0.54  & 0.54  & 0.54  & 0.54  & 0.54  & 0.54 \\
                  & Iter   & 28    & 28    & 28    & 28    & 28    & 28 \\
                  & CPU    & 0.02  & 0.03  & 0.11  & 0.49  & 2.70  & 10.50  \\[0.2cm]

GSOR             & $\ao$   & 0.425 & 0.425 & 0.425 & 0.425 & 0.425 & 0.425 \\
                 & Iter    & 25    & 25    & 25    & 25    & 25    & 25 \\
                 & CPU     & 0.02  & 0.03  & 0.09  & 0.37  & 1.98  & 18.06 \\ \hline
 \end{tabular}
\end{table}

\section{Conclusion}\label{SEC5}

We have presented a two-step iteration method, called TSCSP iteration method, which is a two-step version of the Scale-Splitting (SCSP) iteration method for solving $(W+iT)=b$, where $W$ and $T$ are symmetric positive
semidefinite matrices with at least one of them begin positive definite.
We have shown that if the matrices $W$ and $T$ are symmetric positive definite, then the method is convergent.
We have also obtained the optimal value of the involved parameter $\alpha$. Numerical results show that the TSCSP iteration method often outperforms the SCSP,  the modified HSS (MHSS), the preconditioned MHSS (PMHSS) and  the generalized SOR (GSOR) iteration methods form the number of iterations and the CPU time point of view.

\section*{Acknowledgements}

I would like to thank the three anonymous reviewers and  Prof. M.K. Ng for their constructive comments and suggestions. The work of the author is partially supported by University of Guilan.

\enddocument